\documentclass{article}
\usepackage{hyperref}
\usepackage{ifpdf}
\ifpdf
\usepackage{graphicx}
\usepackage{epstopdf}
\else
\usepackage{graphicx}
   \RequirePackage{breakurl}[1.40]
   \if@nohypdvips
   \else
   \usepackage{hypdvips}
   \fi
\fi
\usepackage{lipsum}
\usepackage{amsfonts}
\usepackage{epstopdf}
\usepackage{algorithmic}
\ifpdf
  \DeclareGraphicsExtensions{.eps,.pdf,.png,.jpg}
\else
  \DeclareGraphicsExtensions{.eps}
\fi

\usepackage{subfigure}
\usepackage{color}
\usepackage{bm}
\usepackage{amssymb}
\usepackage{amsmath}
\usepackage{amsopn}
\def\Xint#1{\mathchoice
   {\XXint\displaystyle\textstyle{#1}}%
   {\XXint\textstyle\scriptstyle{#1}}%
   {\XXint\scriptstyle\scriptscriptstyle{#1}}%
   {\XXint\scriptscriptstyle\scriptscriptstyle{#1}}%
   \!\int}
\def\XXint#1#2#3{{\setbox0=\hbox{$#1{#2#3}{\int}$}
     \vcenter{\hbox{$#2#3$}}\kern-.5\wd0}}
\def\ddashint{\Xint=}

\title{Boundary integral equations for calculating complex eigenvalues of transmission problems\thanks{This work was funded by the Japan Society for the Promotion of Science (JSPS KAKENHI Number 14J03491).}}

\author{
  Ryota Misawa\thanks{Department of Applied Analysis and Complex Dynamical Systems, Kyoto University ({misawa.ryota.27u@kyoto-u.jp, niino@acs.i.kyoto-u.ac.jp, nchml@i.kyoto-u.ac.jp})}
  \and
  Kazuki Niino\footnotemark[2]
  \and
  Naoshi Nishimura\footnotemark[2]
}

\begin{document}

\maketitle

\begin{abstract}
Resonance frequencies are complex eigenvalues at which the homogeneous
transmission problems have non-trivial solutions. These frequencies
are of interest because they affect the behavior of the solutions even
when the frequency is real. The resonance frequencies are related to
problems for infinite domains which can be solved efficiently with the
Boundary Integral Equation Method (BIEM). We thus consider a numerical
method of determining resonance frequencies with fast BIEM and the
Sakurai-Sugiura projection Method (SSM).  However, BIEM may have
fictitious eigenvalues even when one uses M\"uller or PMCHWT
formulations which are known to be resonance free when the frequency
is real valued.

In this paper, we propose new BIEs for transmission problems with
which one can distinguish true and fictitious eigenvalues easily.
Specifically, we consider waveguide problems for the Helmholtz
equation in 2D and standard scattering problems for Maxwell's
equations in 3D.  We verify numerically that the proposed BIEs can
separate the fictitious eigenvalues from the true ones in these
problems.  We show that the obtained true complex eigenvalues are
related to the behavior of the solution significantly.  We also show
that the fictitious eigenvalues may affect the accuracy of BIE
solutions in standard boundary value problems even when the frequency is real.

\end{abstract}

 {key words}:
   resonance, transmission problems, boundary integral equations, eigenvalue problems

\section{Introduction}
There is no doubt that eigenvalue problems are of great importance in
science and engineering. Determining eigenfrequencies of finite
bodies, for example, is one of fundamental issues in the study of
vibrations of structures. In many cases one is interested in real valued
eigenvalues because the physics of the problems requires real valued
quantities as in the case of eigenfrequencies. However, determining
complex eigenvalues is of interest in some applications
even when the quantity in question is real valued physically.  In
waveguides for example, there exist anomalous frequencies near which
the behavior of the solution changes suddenly. Some of such anomalies
are known to be related to the existence of resonance frequencies,
which are complex eigenfrequencies at which there exist non-trivial
solutions to the homogeneous boundary value problems for waveguides.
Interestingly, these resonance frequencies include real valued
eigenvalues called trapped modes which are typically excited near
inclusions (See the review paper by Linton and McIver \cite{Linton_McIver} for efforts
to determine these eigenvalues). It is also known that open systems
may have complex eigenvalues. These complex eigenvalues are called
leaky modes since they radiate energy to infinity. The leaky modes are
also of interest because they are known to affect the behavior of the
solutions considerably.  Actually, numerical examples in this paper
will provide further evidence of the relevance of the leaky modes to
physical phenomena.

In the present paper, we call these problems of finding complex
eigenvalues associated with Boundary Value Problems (BVPs) ``resonance problems''. Specifically, we consider complex eigenfrequencies
(resonance frequencies) for the free space or waveguides which contain
finite sized inclusions.  The complex eigenvalues for our resonance
problems may include real valued ones as in the case of the trapped
modes. Exact solutions of resonance problems are available only
in simple cases and we usually have to resort to numerical or approximate methods to
solve them.  Many efforts have been devoted to numerical and approximate solvers for
resonance problems, e.g., \cite{Hein,Duan,Kim_Pasciak,Hein2010,Gopalakrishnan} to mention just a few. 
The present
authors have been interested in solving resonance problems with
Boundary Integral Equation Method (BIEM) formulated with Green's
functions \cite{Misawaetal}. This method solves resonance problems by finding
frequencies at which the discretized homogeneous BIEs have non-trivial
solutions. BIEM is considered to be advantageous in our resonance
problems because many complex eigenvalue problems are associated with
wave problems for infinite domains where one has to deal with
radiation conditions. %
With BIEM, one does not need special tools such
as PML \cite{Kim_Pasciak} to deal with radiation conditions and the discretization is
required only on the boundary of the inclusions. 
However, BIEM in eigenvalue problem has the following three difficulties: 
\begin{enumerate}
\item The computational cost of BIEM is as large as $O(N^2)$ ($N$: Degrees Of Freedom (DOF)). 

\item One has to solve non-linear eigenvalue problems. 

\item One may obtain non-physical eigenvalues. 
\end{enumerate}
The first difficulty is now solved with the development of the fast
BIEMs such as Fast Multipole Method (FMM \cite{Greengard_Rokhlin}), fast direct solvers
(e.g., \cite{Martinsson_Rokhlin,Bebendorf}), etc. The second difficulty is also much alleviated
with the development of good solvers for non-linear eigenvalue
problems such as the Sakurai-Sugiura projection Method (SSM, \cite{Sakurai_Sugiura,Ikegami_Sakurai,Asakura}). The SSM is a non-iterative algorithm which determines eigenvalues
within a given contour $\gamma$ in the complex plane using contour integrals
defined on $\gamma$, as does another well-known eigensolver FEAST \cite{FEAST}. As a matter of fact, the authors have
developed an FMM for Neumann waveguide problems for the Helmholtz
equation in 2D, and solved resonance problems successfully in \cite{Misawaetal}
with the help of the SSM. The third difficulty is caused by the
difference between the eigenvalues of the BVP and those of the
BIE which may include eigenvalues irrelevant to the original
BVP. These spurious eigenvalues of the BIE are called ``fictitious eigenvalues''. Many efforts have been devoted to the development of
BIEs which are free from real valued fictitious eigenvalues such as combined integral equations \cite{Colton_Kress}, the Burton-Miller equation \cite{Burton_Miller}, the PMCHWT (Poggio-Miller-Chang-Harrington-Wu-Tsai) \cite{Chew} and M\"uller \cite{Muller} formulations.
However, few studies have focused on fictitious eigenvalue issues
in complex eigenvalue problems.

This paper discusses a way of dealing with the third difficulty in
transmission problems. We show that a small modification of the BIE
enables us to clearly distinguish between the true eigenvalues and the
fictitious ones. 
Although we have already commented 
on this modification briefly in our previous paper \cite{Misawaetal} and in conference proceedings \cite{ictam},
the full details of this approach appear for the first time in this paper.
We consider transmission problems for both waveguides for the
Helmholtz equation in 2D and standard Maxwell's equations in 3D. As a
matter of fact, related problems in the exterior Neumann problems for
the Helmholtz equation in 3D have been discussed in \cite{Steinbach_Unger2015} which uses the combined integral equations, etc.
To the best of our knowledge, however, remedies for the
transmission problems have not been developed yet.

This paper is organized as follows: In section \ref{sec:formulation}, we formulate the
transmission resonance problems and the corresponding coupled BIEs for
the waveguide problems for the Helmholtz equation in 2D and the
standard scattering problems for the Maxwell equations in 3D. We
consider both M\"uller's and PMCHWT formulations for these 2 cases. We
then identify fictitious eigenvalues for the integral equations
considered and propose BIEs which can distinguish true and fictitious
eigenvalues clearly.  We show, in particular, that eigenvalues of
the M\"uller and PMCHWT formulations for Maxwell's equations are
identical including fictitious ones.  In Section \ref{sec:numericalexamples}, we present some
numerical examples which prove the effectiveness of the proposed
method.  We also show that true complex eigenvalues affect the
behavior of the solution  and that fictitious complex eigenvalues may
deteriorate the accuracy of the solutions for the ordinary boundary
value problems with real frequencies. This paper ends with a few
concluding remarks and future plans.

\section{Formulation}
\label{sec:formulation}
In this section, we formulate transmission resonance problems and derive BIEs discussed in this study.
We consider one single scatterer for simplicity, although the results in the following discussions hold for multiple scatterer cases as well.

\subsection[Transmission resonance problems for waveguides for the Helmholtz equation]{Transmission resonance problems for waveguides for  the Helmholtz equation}
\label{sec:waveguideBVP}
We first consider elastic waves governed by the Helmholtz equation in 2D. Let $P$ be an infinite strip given by $P=[-1/2,1/2]\times \mathbb{R}$.
Also, let $\Omega=\Omega_2\subset P$ be a finite sized scatterer, $\partial \Omega=\partial \Omega_2$ be its boundary and $\Omega_1=P\setminus \overline{\Omega_2}$.
We consider the following homogeneous transmission problem: find $u$ which satisfies the Helmholtz equation
\begin{displaymath}
 \Delta u +k_{\nu}^2u=0 \quad \mbox{in } \Omega_\nu \ (\nu=1,2),
\end{displaymath}
boundary conditions
\begin{equation}\label{eq:bc}
u^+=u^-\, (=u), \quad S_1\frac{\partial u^+}{\partial n}=S_2\frac{\partial u^-}{\partial n}\,(=q)\quad\mbox{on }\partial \Omega, 
\end{equation}
and the homogeneous Neumann boundary condition on the sides of the strip:
\begin{displaymath}
\frac{\partial u}{\partial x_1}=0 \quad\mbox{on } x_1=\pm\frac{1}{2}.
\end{displaymath}
We impose the radiation condition which requires that $u$ is written as follows (e.g. \cite{Linton_Evans}):
\begin{align}
  u(x)\approx\sum_{l\ge 0} C_l^{\pm} \cos l\pi\left(x_1+\frac{1}{2}\right)e^{-\sqrt{(l\pi)^2-k_1^2}|x_2|} \ \mbox{as} \ x_2\rightarrow \pm\infty\label{radiation},
\end{align}
where $\rho_\nu$, $S_\nu$  and $k_\nu=\omega\sqrt{\rho_\nu/S_\nu}$ are the density, shear modulus (real numbers)  and the wavenumber in $\Omega_\nu$ ($\nu=1,2$), respectively. The frequency $\omega$ is allowed to take a complex value. Also, the superscript  $+$ ($-$) in \eqref{eq:bc} stands for the trace to $\partial \Omega$ from $\Omega_1$ ($\Omega_2$), $\partial/\partial n$ for the normal derivative and $\bm n$ for the unit normal vector on $\partial \Omega$ directed towards $\Omega_1$, respectively.
We make the square root $\sqrt{(l\pi)^2-k_1^2}$ in (\ref{radiation}) single valued as a function of $k_1$  by taking the branch
which is analytic in the complex plane
cut along $(-l\pi,-l\pi-i\infty)$ and $(l\pi,l\pi-i\infty)$ and
approaches $-ik_1$ in the upper plane as $|k_1|\to \infty$.
This definition of the square root ensures that the summands in (\ref{radiation}) decay as $l\rightarrow\infty$ and allows the analytic continuation of the radiation condition to a complex $\omega$.
We also note that the behavior of the solution in the far field can be described by a finite sum with $l\le \Re k_1/\pi$ in \eqref{radiation} since terms with $l>\Re k_1/\pi$ go to zero as $x_2\rightarrow \pm \infty$.

In the following, we call the above homogeneous problem the ``waveguide problem'' or the ``original BVP''.
The transmission resonance problem determines frequencies at which the waveguide problem has non-trivial solutions.
We call such frequencies ``true eigenvalues''. We note that the true eigenvalues of the waveguide problem have non-positive imaginary parts.
\subsection{BIEs for the Helmholtz equation}
We now formulate BIEs for the transmission resonance problem.
We define $U^\nu$ as follows:
 \begin{align}
U^\nu(x)=(-1)^\nu\left(\frac{1}{S_\nu} \int_{\partial \Omega} G^\nu(x,y) q(y)ds_y-\int_{\partial \Omega} \frac{\partial G^\nu(x,y)}{\partial n_y} u(y) ds_y\right), \ \nu=1,2\label{pot}
\end{align}
where $G^1$ stands for Green's function for the waveguide with the wavenumber $k_1$  and $G^2$ for the fundamental solution with the wavenumber $k_2$, respectively:
\begin{align}
  & G^1(x,y)=\sum_{l=0}^{\infty} f_l \frac{e^{-\sqrt{(l\pi)^2-k_1^2}|x_2-y_2|}}{\sqrt{(l\pi)^2-k_1^2}}\cos{l\pi\left(x_1+\frac{1}{2}\right)}\cos{l\pi\left(y_1+\frac{1}{2}\right)}\label{green}\\
 &G^2(x,y)=\frac{i}{4} H_0^{(1)}(k_2|\bm x-\bm y|)\label{fund}.
\end{align}
In (\ref{green}),  $f_l$ is 1 for $l\ne 0$ and  1/2 for $l=0$ and $H_{n}^{(\iota)}$ stands for the Hankel function of the $\iota$-th kind and $n$-th order, respectively.   Also, bold letters $\bm x$, $\bm y$, etc.\ stand for the position vectors of the points $x$, $y$, etc.

It is well-known that $U^\nu$ ($\nu=1,2$) give potential representations of the solution $u$ of the transmission problem in \ref{sec:waveguideBVP} in $\Omega_\nu$ with the boundary traces of the solution as the layer-potential densities, iff $U^1$ ($U^2$) vanishes in $\Omega_2$ ($\mathbb{R}^2\setminus \overline{\Omega_2}$). This condition gives 
\begin{align}
  U^{1-}=U^{2+}=\frac{\partial U^{1-}}{\partial n}=\frac{\partial U^{2+}}{\partial n}=0 \ \mbox{on } \partial \Omega.
\label{poteq0}
\end{align}
The M\"uller formulation  of BIE for the Helmholtz equation in 2D given by
\begin{equation}
\begin{split}
  &    \frac{S_1+S_2}{2} u-\int_{\partial \Omega} \left(S_1 \frac{\partial G^1}{\partial n_y} -S_2\frac{\partial G^2}{\partial n_y}\right) u \, ds_y+\int_{\partial \Omega} \left(G^1-G^2\right) q \, ds_y=0\\
  & \frac{S_1+S_2}{2S_1S_2} q-\int_{\partial \Omega} \left(\frac{\partial^2 G^1}{\partial n_x\partial n_y} -\frac{\partial^2 G^2}{\partial n_x\partial n_y}\right) u \, ds_y+\int_{\partial \Omega} \left(\frac{1}{S_1} \frac{\partial G^1}{\partial n_x} -\frac{1}{S_2}\frac{\partial G^2}{\partial n_x}\right) q \, ds_y=0
\end{split}
\label{muller_explicit}
\end{equation}
is derived from (\ref{poteq0}) with the help of the conditions
\begin{align}
 -S_1U^{1-}=S_2U^{2+}, \ 
  -\frac{\partial U^{1-}}{\partial n}=\frac{\partial U^{2+}}{\partial n}\label{muller} \quad \mbox{on } \partial \Omega.
\end{align}
We also consider the following PMCHWT  formulation:
\begin{equation}
\begin{split}
  &    \int_{\partial \Omega} \left(\frac{\partial G^1}{\partial n_y}+\frac{\partial G^2}{\partial n_y}\right) u \, ds_y-\int_{\partial \Omega} \left(\frac{1}{S_1}G^1+\frac{1}{S_2}G^2\right) q \, ds_y=0\\
  & \ddashint_{\partial \Omega} \left(S_1\frac{\partial^2 G^1}{\partial n_x\partial n_y} +S_2\frac{\partial^2 G^2}{\partial n_x\partial n_y}\right) u \, ds_y-\int_{\partial \Omega} \left(\frac{\partial G^1}{\partial n_x}+\frac{\partial G^2}{\partial n_x}\right) q \, ds_y=0
\label{pmchwt_explicit}
\end{split}
\end{equation}
obtained similarly from the conditions
\begin{align}
U^{1-}=U^{2+}, \ 
   S_1\frac{\partial U^{1-}}{\partial n}= S_2\frac{\partial U^{2+}}{\partial n}\label{pmchwt} \quad \mbox{on } \partial \Omega,
\end{align}
where $\ddashint$ stands for the finite part of a divergent integral.

One may want to solve the transmission resonance problems by finding complex $\omega$'s such that these integral equations have non-trivial solutions, i.e.,  the non-linear eigenvalues for BIEs in  (\ref{muller_explicit}) or (\ref{pmchwt_explicit}).
This issue is further discussed in the next section.

\subsection{Fictitious eigenvalues}
We have reduced the transmission resonance \\ 
problem to a non-linear eigenvalue problem for BIEs in  (\ref{muller_explicit}) or (\ref{pmchwt_explicit}).  However, as we shall see, these equations may pick up fictitious eigenvalues which we now characterize.

  We discuss in detail the M\"uller formulation because the PMCHWT case can be treated in a similar manner.
We note that the following statement holds (a similar statement is given in \cite{Shipman2003}): A frequency $\omega$ at which the BIEs in (\ref{muller_explicit}) have non-trivial solutions corresponds either to the eigenvalue of the original waveguide problem or  that of the following transmission resonance problem for $v$ in the free space $\mathbb{R}^2$ in which the governing equations in  $\Omega_1$ and $\Omega_2$ are interchanged; we refer to \cite{Kim_Pasciak} for the (3-D) explicit form of the radiation condition for a complex wavenumber:
\begin{align}
 &  \Delta v+k_2^2v=0 \ \mbox{in} \ \mathbb{R}^2\setminus \overline{\Omega_2}, \quad \Delta v+k_1^2v=0 \ \mbox{in} \ {\Omega_2} ,\label{fictbvp_helm} \\
 &v^{-}=v^{+}, \quad S_2\frac{\partial v^-}{\partial n}=S_1\frac{\partial v^+}{\partial n} \quad \mbox{on } \partial \Omega \label{fictbvp_bc}\\
 &\mbox{outgoing radiation condition with } k_2 \mbox{ in } \mathbb{R}^2\setminus\overline{\Omega_2}.
 \label{fictbvp_rad}
\end{align}
To see this, we define a function $v$ using a set of non-trivial solutions $(u,q)$ of the BIEs in (\ref{muller_explicit}) as:
\begin{align}
  v(x)=\left\{
    \begin{array}{c}
     -\frac{1}{S_1}U^2(x)  \quad x \in \mathbb{R}^2\setminus \overline{\Omega_2}\\
     \frac{1}{S_2}U^1(x) \quad x \in \Omega_2
    \end{array}
    \right.
    \label{sppot}.
\end{align}
If $v\equiv 0$ in $\mathbb{R}^2\setminus \partial \Omega$, the function $w$ defined by:
\begin{align}
  w(x)=\left\{
    \begin{array}{c}
      U^1(x)  \quad x \in \Omega_1\\
      U^2(x) \quad x \in \Omega_2
    \end{array}
    \right.
    \label{sppot2}.
\end{align} 
gives a non-trivial solution of the waveguide problem, as we have noted.
We thus see that the $\omega$ is an eigenvalue of the waveguide problem.
If $v \not \equiv 0$ identically in $\mathbb{R}^2\setminus\partial \Omega$, $v$ is a non-trivial solution of the  BVP in (\ref{fictbvp_helm})-(\ref{fictbvp_rad}) as we can see from (\ref{muller}). 
We call this BVP the ``fictitious BVP''.

A similar discussion shows that the PMCHWT formulation given in (\ref{pmchwt_explicit}) may have, in addition to true eigenvalues, fictitious eigenvalues which correspond to eigenvalues of the following free space transmission problem:
\begin{align}
 &  \Delta v+k_2^2v=0 \ \mbox{in} \ \mathbb{R}^2\setminus \overline{\Omega_2}, \quad \Delta v+k_1^2v=0 \ \mbox{in} \ {\Omega_2} ,\label{fictbvp_helm2} \\
 &v^{-}=v^{+}, \quad S_1\frac{\partial v^-}{\partial n}=S_2\frac{\partial v^+}{\partial n}\ \mbox{on }\partial \Omega\label{fictbvp_bc2}\\
 &\mbox{outgoing radiation condition with } k_2 \mbox{ in } \mathbb{R}^2\setminus\overline{\Omega_2}.
  \label{fictbvp_rad2}
\end{align}
One may say that fictitious eigenvalues exist because the BIEs in (\ref{muller_explicit}) and (\ref{pmchwt_explicit}) cannot distinguish the original BVP and the fictitious BVP since these BIEs lose the information of the domains while we take their limits to the boundary.

We thus see that an eigenvalue obtained with the BIEM with (\ref{muller_explicit}) may be an eigenvalue of the free space transmission problem in (\ref{fictbvp_helm})-(\ref{fictbvp_rad}), which has nothing to do with the original BVP and, hence, is fictitious.   
 These fictitious eigenvalues have negative imaginary parts because the homogeneous transmission problem has only the trivial solution when $\Im \omega\ge 0$  \cite{Kress}.
Hence, the BIE in (\ref{muller_explicit}) is free from fictitious eigenvalues as long as one considers real frequencies. When one deals with leaky modes, however, it is  hard to tell
 whether an eigenvalue obtained with the BIEMs in (\ref{muller_explicit}) is a true eigenvalue or not, because both true and fictitious eigenvalues appear in the lower complex plane. Similar conclusions apply to (\ref{pmchwt_explicit}) as well.

\subsection{New BIEs}
One can resolve the above problem simply by  using  the incoming fundamental solution given by:
\begin{align}
 -\frac{i}{4} H_0^{(2)}(k_2|\bm x-\bm y|)\label{incomingfund}
\end{align}
for $G^{2}$ in (\ref{muller_explicit}), instead of the outgoing one in (\ref{fund}). This remedy keeps the true eigenvalues unchanged, while the corresponding fictitious transmission problem (\ref{fictbvp_helm})-(\ref{fictbvp_rad}) is now replaced by the one with (\ref{fictbvp_helm}), (\ref{fictbvp_bc}) and the incoming radiation condition with $k_2$ in $\mathbb{R}^2\setminus\overline{\Omega_2}$.  The corresponding (fictitious) eigenvalues have positive imaginary parts because of the incoming radiation condition. In fact, we see that $(\overline{v},\bar{\omega})$ is 
an eigenpair of the problem defined by (\ref{fictbvp_helm}), (\ref{fictbvp_bc}) and the incoming radiation condition  if $(v,\omega)$ is an eigenpair  of  (\ref{fictbvp_helm})-(\ref{fictbvp_rad}).
Therefore, we can distinguish the fictitious eigenvalues from the true eigenvalues with this change of the formulation. The same method can be applied to (\ref{pmchwt_explicit}) as well.

\subsection{Transmission resonance problems for the Maxwell equations}
We next consider transmission problems for Maxwell's equations in 3D free space. Let $\Omega_2$ be a finite scatterer, $S=\partial \Omega_2$ be its boundary and $\Omega_1=\mathbb{R}^3\setminus \overline{\Omega_2}$.  The transmission resonance problem for the Maxwell equations is stated as follows: find $\bm E$ which satisfies the Maxwell equations:
\begin{displaymath}
  \nabla\times \left(\nabla\times \bm E\right)-k_\nu^2\bm E=0 \quad \mbox{in } \Omega_\nu  \ (\nu=1,2)
\end{displaymath}
boundary conditions
\begin{equation}
\begin{split}
 & \bm E^+\times \bm n=\bm E^-\times \bm n\ (=\bm m),\\ 
&\bm n\times \frac{1}{i\omega} \left(\frac{1}{\mu_1}\nabla\times \bm E^{+} \right)=\bm n\times \frac{1}{i\omega} \left(\frac{1}{\mu_2}\nabla\times \bm E^{-} \right)\ (=\bm j)
 \quad \mbox{on } S
\end{split}
\label{bc_maxwell}
\end{equation}
and the outgoing radiation condition with $k_1$ in $\Omega_1$ given by:
\begin{equation}
  \begin{split}
  \bm E=\sum_{n=1}^{\infty} \sum_{m=-n}^{n} \alpha_n^m \bm \Phi_n^m(\bm x)+\beta_n^m \nabla\times \bm \Phi_n^m(\bm x), \ \bm \Phi_n^m(\bm x)=\nabla\times \left\{ \bm x h_n^{(1)}(k_1r)Y_n^m(\theta,\phi)\right\}\\
  \mbox{for } |\bm x|>R
\label{radiation_maxwell}
  \end{split}
\end{equation}
where superscript  $+$ ($-$) stands for the trace to $S$ from $\Omega_1$ ($\Omega_2$), $\epsilon_\nu,\ \mu_\nu$ are the permitivity and permeability (real numbers) and $k_\nu=\omega\sqrt{\epsilon_\nu\mu_\nu}$ in $\Omega_\nu$, respectively.  Also, $\alpha_n^m$ and $\beta_n^m$ are numbers, $(r,\theta,\phi)$ is the spherical coordinate of $\bm x$, $h_n^{(1)}$ is the spherical Hankel function of the 1st kind, $Y_n^m$ is the spherical harmonics and $R>0$ is a constant such that $|\bm x|<R$ holds for $\forall\bm x\in \overline{\Omega_2}$.
The above expression in (\ref{radiation_maxwell})  allows the analytic continuation of the Silver--M\"uller radiation condition to a complex $k_1$.
\subsection{BIEs for the Maxwell equations}
We next consider the BIEs for the Maxwell equations.
We introduce the following potential representation of $\bm E$ via the surface magnetic current $\bm m$ and the electronic current $\bm j$ (see (\ref{bc_maxwell})):
\begin{equation}
\begin{split}
  &\tilde{E}^{\nu}_i(x)\\
  &=(-1)^\nu \int_{S} \left\{e_{ijk} \frac{\partial \Gamma^{\nu}(x,y)}{\partial x_j} m_k(y)
  -i\omega \mu_{\nu} \left(\delta_{ip}+\frac{1}{k^2_{\nu}}\frac{\partial}{\partial y_i}\frac{\partial}{\partial y_p}\right)\Gamma^{\nu}(x,y)j_p(y) \right\}\, dS_y\label{pot_maxwell}
\end{split}
\end{equation}
where $\Gamma^\nu(x,y)$ stands for the fundamental solution of the Helmholtz equation in 3D:
\begin{align}
  \Gamma^\nu(x,y)=\frac{e^{ik_{\nu}|\bm x-\bm y|}}{4\pi |\bm x-\bm y|} \label{fund3d}.
\end{align}
These representations give the solution of the original BVP with the boundary traces of the solution as the densities $\bm m$ and $\bm j$ iff the following equations hold:
\begin{eqnarray}
\tilde{\bm{E}}^1=0 \quad \mbox{in } \Omega_2, \ 
\tilde{\bm{E}}^2=0 \quad \mbox{in } \Omega_1.
\label{poteq0_maxwell}
\end{eqnarray}

The integral equations for the M\"uller formulation can be written  as follows \cite{Muller}:
\begin{equation}
\begin{split}
 & \frac{\epsilon_1+\epsilon_2}{2} \bm m-\bm n\times \int_{S} \left(
\epsilon_1\nabla_x \Gamma^1-\epsilon_2\nabla_x \Gamma^2\right)\times \bm m \, dS_y\\
&+
\frac{i}{\omega} \bm n\times \int_S \left(k_1^2\Gamma^1-k_2^2\Gamma^2\right)\bm j\, dS_y+
\frac{i}{\omega}\bm n\times \int_S \left(\nabla_x\nabla_x \Gamma^1-\nabla_x \nabla_x \Gamma^2\right)\cdot \bm j\, dS_y=0,\\
 & \frac{\mu_1+\mu_2}{2} \bm j-\bm n\times \int_{S} \left(
\mu_1\nabla_x \Gamma^1-\mu_2\nabla_x \Gamma^2\right)\times \bm j \, dS_y\\
&-
\frac{i}{\omega} \bm n\times \int_S \left(k_1^2\Gamma^1-k_2^2\Gamma^2\right)\bm m\, dS_y-
\frac{i}{\omega}\bm n\times \int_S \left(\nabla_x\nabla_x \Gamma^1-\nabla_x \nabla_x \Gamma^2\right)\cdot \bm m\, dS_y=0. 
\end{split}
\label{muller_maxwell_explicit}
\end{equation}
The above M\"uller formulation for the the Maxwell equations is derived from (\ref{poteq0_maxwell})  via the potential in (\ref{pot_maxwell}) as follows:
\begin{align}
  \begin{array}{c}
    -\epsilon_1 \tilde{\bm{E}}^{1-}\times \bm n=\epsilon_2 \tilde{\bm{E}}^{2+}\times \bm n\\
-\frac{1}{i\omega} \bm n\times \left(\nabla\times {\tilde{\bm E}}^{1-}\right)=\frac{1}{i\omega} \bm n\times \left(\nabla\times {\tilde{\bm E}}^{2+}\right)
\end{array}
\quad
\mbox{on }
S.
\label{muller_maxwell}
\end{align}

The PMCHWT formulation  for the  Maxwell equations can be written as follows \cite{Chew}:
\begin{equation}
\begin{split}
 &-\bm n\times \int_{S} \left(
\nabla_x \Gamma^1+\nabla_x \Gamma^2\right)\times \bm m \, dS_y
+
{i\omega} \bm n\times \int_S \left(\mu_1\Gamma^1+\mu_2\Gamma^2\right)\bm j\, dS_y\\
&+
\frac{i}{\omega}\bm n\times \ddashint_S \left(\frac{1}{\epsilon_1}\nabla_x\nabla_x \Gamma^1+\frac{1}{\epsilon_2}\nabla_x \nabla_x \Gamma^2\right)\cdot \bm j\, dS_y=0,\\
 & -\bm n\times \int_{S} \left(
\nabla_x \Gamma^1+\nabla_x \Gamma^2\right)\times \bm j \, dS_y
-
{i\omega} \bm n\times \int_S \left(\epsilon_1\Gamma^1+\epsilon_2\Gamma^2\right)\bm m\, dS_y\\
&-
\frac{i}{\omega}\bm n\times \ddashint_S \left(\frac{1}{\mu_1}\nabla_x\nabla_x \Gamma^1+\frac{1}{\mu_2}\nabla_x \nabla_x \Gamma^2\right)\cdot \bm m\, dS_y=0,
\end{split}
\label{pmchwt_maxwell_explicit}
\end{equation}
which we obtain from (\ref{pot_maxwell}) and (\ref{poteq0_maxwell})  using 
\begin{align}
\begin{array}{c}
    \tilde{\bm{E}}^{1-}\times \bm n=\tilde{\bm{E}}^{2+}\times \bm n\\
\frac{1}{i\omega} \bm n\times \left(\frac{1}{\mu_1}\nabla\times \tilde{\bm E}^{1-}\right)=\frac{1}{i\omega} \bm n\times \left(\frac{1}{\mu_2}\nabla\times \tilde{\bm E}^{2+}\right)
\end{array}
\quad
\mbox{on }
S.
\label{pmchwt_maxwell}
\end{align}
 The above homogeneous BIEs in (\ref{muller_maxwell_explicit}) and (\ref{pmchwt_maxwell_explicit}) have both true eigenvalues and fictitious ones as in the Helmholtz case. The fictitious eigenvalues for the M\"uller formulation are the eigenvalues $\omega$ of the following BVP for $\bm H$:
\begin{align}
  &\nabla\times \left(\nabla\times \bm H\right)-k_\nu^2\bm H=0
  \quad \mbox{in} \quad \Omega_{\nu'}  \ (\nu\ne\nu', \ \nu=1,2)
  \label{maxwelleq_fict_h}\\
  &\bm H^{-}\times \bm n=\bm H^{+} \times \bm n, \quad
  \bm n \times \left(\frac{1}{\epsilon_1}\nabla\times \bm H^{-}\right)=\bm n \times \left(\frac{1}{\epsilon_2}\nabla\times\bm H^{+}\right) \quad \mbox{on } S\label{bc_m}\\
 &\mbox{outgoing radiation condition for }\bm H\mbox{ with } k_2 \mbox{ in } {\Omega_1}
\label{radiation_fict_h}
\end{align}
where $\bm H$ is related to $\tilde{\bm E}$ in (\ref{muller_maxwell}) via
\begin{align}
  \bm H=\left\{ 
  \begin{array}{c}
    -\epsilon_1 {\tilde{\bm E}}^1 \ \mbox{in } \Omega_2\\
    \epsilon_2 {\tilde{\bm E}}^2 \ \mbox{in } \Omega_1
    \end{array}
    \right. .
  \end{align}
The fictitious boundary value problem for the PMCHWT formulation is given by 
\begin{align}
&  \nabla\times \left(\nabla\times \bm E\right)-k_\nu^2\bm E=0
  \quad \mbox{in} \quad \Omega_{\nu'}  \ (\nu\ne\nu', \ \nu=1,2)\label{maxwelleq_fict_e}\\
 & \bm E^{-}\times \bm n=\bm E^{+} \times \bm n, \,
  \bm n \times \left(\frac{1}{\mu_1} \nabla\times \bm E^{-}\right)=\bm n \times \left(\frac{1}{\mu_2} \nabla\times \bm E^{+}\right)\,\mbox{on } S\label{bc_p}\\
 &\mbox{outgoing radiation condition for }\bm E\mbox{ with } k_2 \mbox{ in } {\Omega_1}.
\label{radiation_fict_e}
\end{align}
where $\bm E={\tilde{\bm E}}^\nu \ \mbox{in }\Omega_{\nu'}$.
As a matter of fact, the  fictitious eigenvalues for the M\"uller  and PMCHWT formulations are the same. Indeed, the fictitious BVPs for the M\"uller and  PMCHWT formulations transform to each other by the following ``changes of variables'':
\begin{align}
  &\bm H=\frac{\nabla \times \bm E}{i\omega \mu_\nu} \ \mbox{in } \Omega_{\nu'}\label{change_of_variables1}
\\
  &\bm E=-\frac{\nabla \times \bm H}{i\omega \epsilon_\nu} \ \mbox{in } \Omega_{\nu'}
\label{change_of_variables2}
\end{align}
Namely, one obtains (\ref{bc_p}) by rewriting (\ref{bc_m})  with (\ref{change_of_variables2}). Also, one obtains (\ref{bc_m}) by using (\ref{change_of_variables1}) in (\ref{bc_p}). 
These changes of variables also keep the Maxwell equations and the radiation conditions unchanged. Hence the eigenvalues of the M\"uller and PMCHWT formulations are identical including fictitious ones.
Incidentally, this conclusion is quite obvious from a physical point of view since the problems defined by (\ref{maxwelleq_fict_h})-(\ref{radiation_fict_h}) and (\ref{maxwelleq_fict_e})-(\ref{radiation_fict_e}) are the same transmission problem formulated in terms of either the magnetic or electric field.

Also in Maxwell's equations with (\ref{muller_maxwell_explicit}) or (\ref{pmchwt_maxwell_explicit}),  we can distinguish true and fictitious eigenvalues by replacing $\Gamma^2$ with the incoming fundamental solution given by:
\begin{align*}
  \frac{e^{-ik_2 |\bm x-\bm y|}}{4\pi|\bm x-\bm y|}.
\end{align*}
A fictitious eigenpair $(\bm H, \omega)$ of (\ref{maxwelleq_fict_h})-(\ref{radiation_fict_h}) are changed to $(\bar{\bm H},\bar{\omega})$ with this new formulation.

\label{sec:maxwell_fict}

\section{Numerical examples}
\label{sec:numericalexamples}
In this section, we present some numerical examples to test the performances of the proposed method. We used Appro GreenBlade 8000 (with Intel Xeon cores) at the Academic Center for Computing and Media Studies of Kyoto University, and FX10 Supercomputer System (with $\mbox{SPARC64}^{\mathrm{TM}}\mbox{IXfx}$ cores) at the Information Technology Center of the University of Tokyo for the computation. The codes are parallelized with OpenMP and MPI.
\subsection{Discretization of BIEs}
We use both M\"uller and PMCHWT formulations (in (\ref{muller_explicit}) and (\ref{pmchwt_explicit}), respectively) for solving the waveguide problems for the Helmholtz equation. 
The BIEs  in  (\ref{muller_explicit}) and (\ref{pmchwt_explicit})  are discretized  with 
piecewise constant boundary elements and the collocation method (the singular parts of the integrals are evaluated analytically, and the remainders are computed with the Gaussian quadrature).
The discretized BIE for (\ref{muller_explicit}) converges fast   in GMRES (Generalized Minimal RESidual method, \cite{Saad}) since the operator is a compact perturbation of a constant.
The PMCHWT formulation in (\ref{pmchwt_explicit}) discretized with collocation is also known to converge fast in GMRES  if the unknowns are ordered in a proper manner (\cite{Niino2d}, \cite{Niino3d}).
The matrix-vector product operation is accelerated with the FMM for waveguide problems proposed in \cite{Misawaetal}.

For the Maxwell equations, we use only the M\"uller formulation in (\ref{muller_maxwell_explicit})  because eigenvalues of the M\"uller and PMCHWT formulations are identical as we have noted in \ref{sec:maxwell_fict}. 
The BIE in \eqref{muller_maxwell_explicit} is discretized with triangular boundary elements and Nystr\"om's method discussed in \cite{Niino_Muller} using  the three point  Gaussian quadrature rule on a triangle.   For the local correction, the contributions of the static parts of the fundamental solution are calculated analytically after interpolating the densities linearly, and those of the remainder are integrated numerically with the Gaussian quadrature.  %
The discretized BIE converges fast with GMRES \cite{Niino_Muller}. The matrix-vector product operation is accelerated with the low frequency FMM following \cite{Niino_Muller}.

\subsection{Sakurai-Sugiura projection method (SSM)}
We briefly describe the SSM \cite{Sakurai_Sugiura,Ikegami_Sakurai,Asakura} for solving non-linear eigenvalue problems for BIEs.  We write the discretized version of the homogeneous BIEs as
\begin{align}
  A(\omega) \bm x=\bm 0 \label{ax}
\end{align}
where $A$ represents the $N\times N$ matrix of the discretized BIEs which depends on complex $\omega$ in a non-linear manner.  
Our eigenvalue problem finds  $\omega$'s at which (\ref{ax}) has non-trivial solutions.
The SSM determines eigenvalues within a  given contour $\gamma$ in the complex plane.

The SSM reduces the non-linear eigenvalue problem for $A(\omega)$ to a generalized eigenvalue problem given by $H_{mL}^{<}\bm x=z H_{mL}\bm x$ for two Hankel matrices defined as  
$H_{mL}=\left[M_{i+j-2}\right]_{i,j=1}^m\in \mathbb{C}^{mL\times mL}$ and $H_{mL}^{<}=\left[M_{i+j-1}\right]_{i,j=1}^m\in \mathbb{C}^{mL\times mL}$, where $M_k$ is a matrix defined as follows: 
\begin{align}
  M_k=P^H\frac{1}{2\pi i}\int_{\gamma} z^k  A^{-1}(z)Q dz, \ k=0,\cdots, 2m-1\label{mu},
\end{align}
and  $P$ and $Q$ are random $N\times L$ matrices (namely, $L$ random vectors of dimension $N$), respectively. The non-linear eigenvalues of $A$ are then obtained as the solutions of the above generalized eigenvalue problem \cite{Asakura}. The integer parameters $m$ and $L$ are set so that $mL$ is larger than the number of eigenvalues within $\gamma$ (our typical choice of $m$ is 12 or 24 in the examples to follow).

When we solve the waveguide problem with Green's function in (\ref{green}), however, we need to take the path of integration $\gamma$  so that it does not touch the branch cuts of Green's function \cite{Misawaetal} given by
\begin{align}
  k_1=p\pi-i\kappa, \ p\in \mathbb{N}, \kappa\ge 0\label{branchcuts}.
\end{align}

We set the contour $\gamma$ of the SSM to be a rectangle in the following numerical examples.  We compute the contour integral in (\ref{mu}) numerically  with the standard Gaussian quadrature applied to 4 integrals on each side of the rectangle.  The inverse of the matrix in (\ref{mu}) is computed approximately with the FMM and GMRES. We set the tolerance (relative error) for GMRES to be $10^{-8}$ in the examples to follow.

\label{ssmethod}

\subsection{Resonances in 2D Helmholtz waveguides}
We first discuss waveguide problems for the Helmholtz equation in 2D defined in \ref{sec:waveguideBVP}. We consider 4 circular scatterers with the radii of $r_0=0.4$ whose centers are $(0,0)$, $(0,-1)$, $(0,-2)$ and $(0,-3)$, respectively.  In this case, we can check if the proposed approach is able to separate fictitious eigenvalues from the true ones since the fictitious eigenvalues can be determined semi-analytically.
In fact,  the fictitious eigenvalues $\omega$ for 
M\"uller's integral equation in (\ref{muller_explicit}), are zeros of  the following expression 
\begin{align}
  -S_2H_n^{(\iota)}(k_2r_0)\frac{d}{dr} J_n(k_1r_0)+S_1\frac{d}{dr}H_n^{(\iota)}(k_2r_0)J_n(k_1r_0)
\label{mie_helmholtz}
\end{align}
where $\iota=1$ for the ordinary method (outgoing) and $\iota=2$ for the proposed method (incoming), respectively.
The fictitious eigenvalues for the PMCHWT case are  calculated similarly.

The following material constants are used in waveguide problems: $\rho_1=1$, $S_1=1$, $\rho_2=0.37$ and $S_2=0.2$. They are supposed to model flint glass inclusions within an iron plate with an appropriate normalization. We discretized each circular boundary with 4000 elements (DOF is 32000), used 32 (64) integration points on each side of $\gamma$ to approximate the integrals in (\ref{mu}) for $\Re \omega<\pi$ ($\Re\omega>\pi$) and set $L=10$ ($L$: number of random vectors used in SSM. See \ref{ssmethod}.).  The number of  integration points on $\gamma$  is chosen so that we can  calculate eigenvalues close to the branch cuts accurately. 
We tested 4 methods shown in Table \ref{tbl:bies}. The methods 1 and 3 are the proposed methods while methods 2 and 4 are standard ones.
\begin{table}
  \caption{BIEs}
  \begin{center}
    \begin{tabular}{|c|c|c|}\hline
      Method No.&kernel function $G^2$ in (\ref{pot})& formulation\\
      \hline
      1& $(-i/4)H_0^{(2)}(k_2|\bm x-\bm y|)$& M\"uller\\
      2& $(i/4)H_0^{(1)}(k_2|\bm x-\bm y|)$& M\"uller\\
      3& $(-i/4)H_0^{(2)}(k_2|\bm x-\bm y|)$& PMCHWT\\
      4& $(i/4)H_0^{(1)}(k_2|\bm x-\bm y|)$& PMCHWT\\
      \hline
    \end{tabular}
  \end{center}
    \label{tbl:bies}
\end{table}

\begin{figure}[]
   \begin{center}
     {\includegraphics[width=\textwidth]{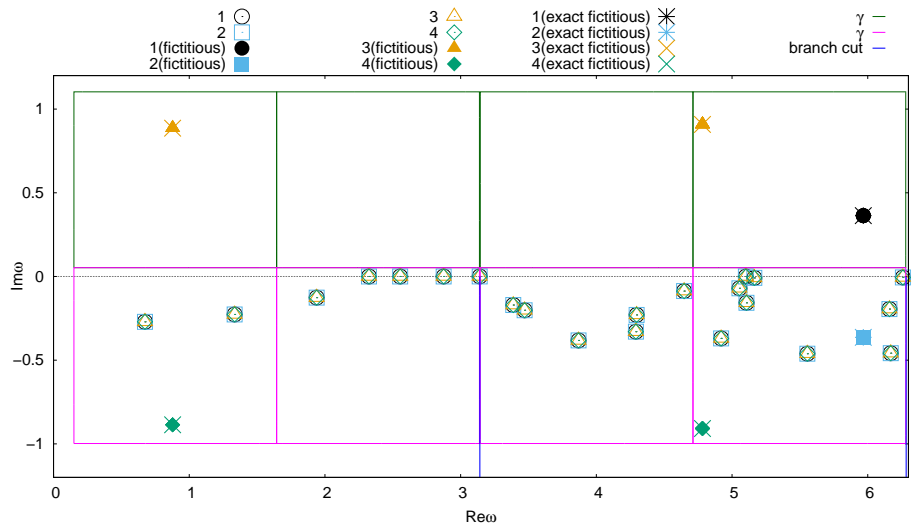}}
     \caption{Eigenvalues obtained with each method. Solid symbols stand for fictitious eigenvalues. The paths $\gamma$'s are taken to be close to the branch cuts in (\ref{branchcuts}) but not to touch them.}
     \label{fig:eig}
   \end{center}
 \end{figure}

\begin{figure}[]
   \begin{center}
     {\includegraphics[width=\textwidth]{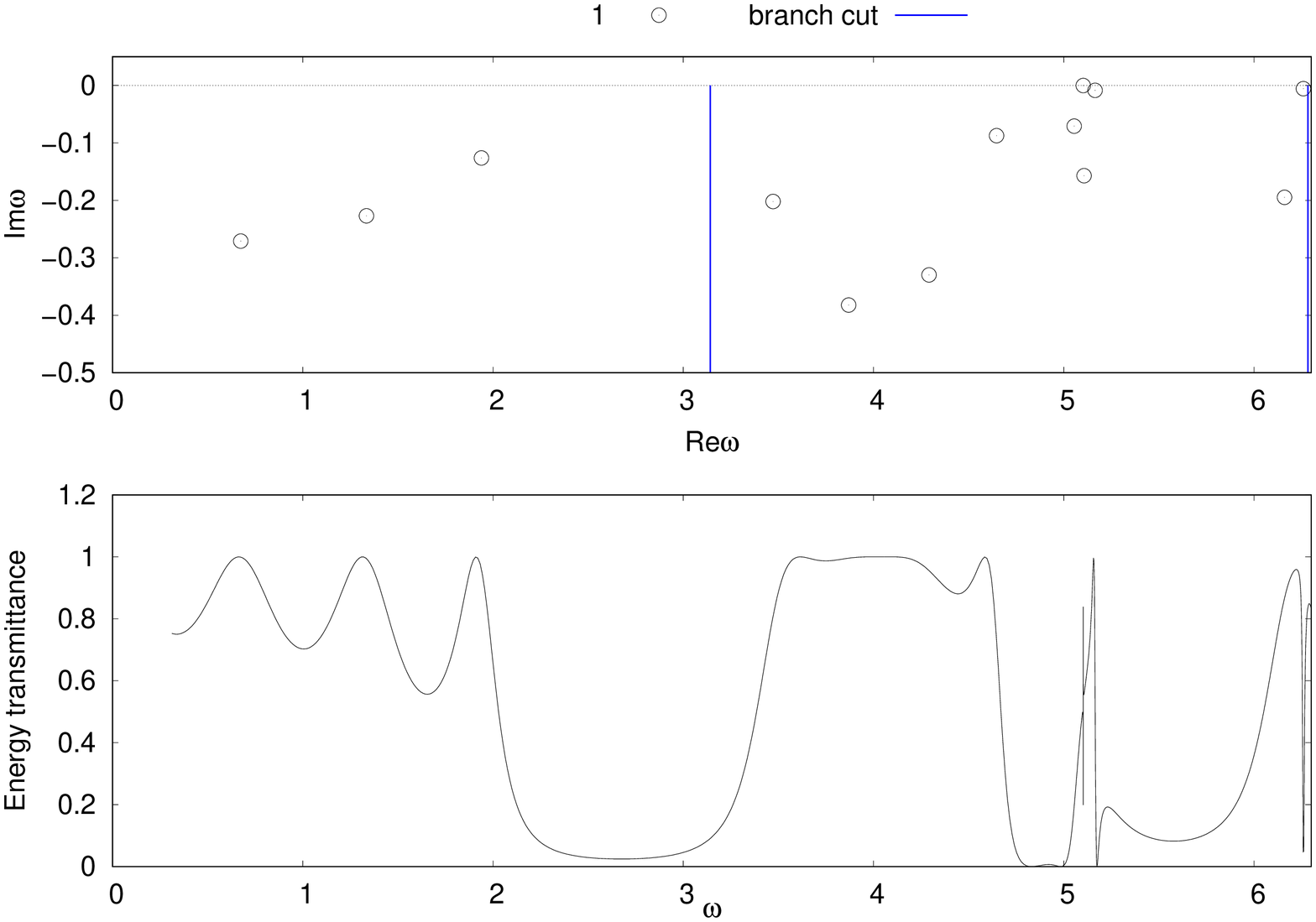}}
     \caption{Eigenvalues with symmetric modes (upper) and energy transmittance for incident plane wave $e^{ik_1x_2}$ (lower).}
     \label{fig:eig_even}
   \end{center}
 \end{figure}

\begin{figure}[]
   \begin{center}
     {\includegraphics[width=\textwidth]{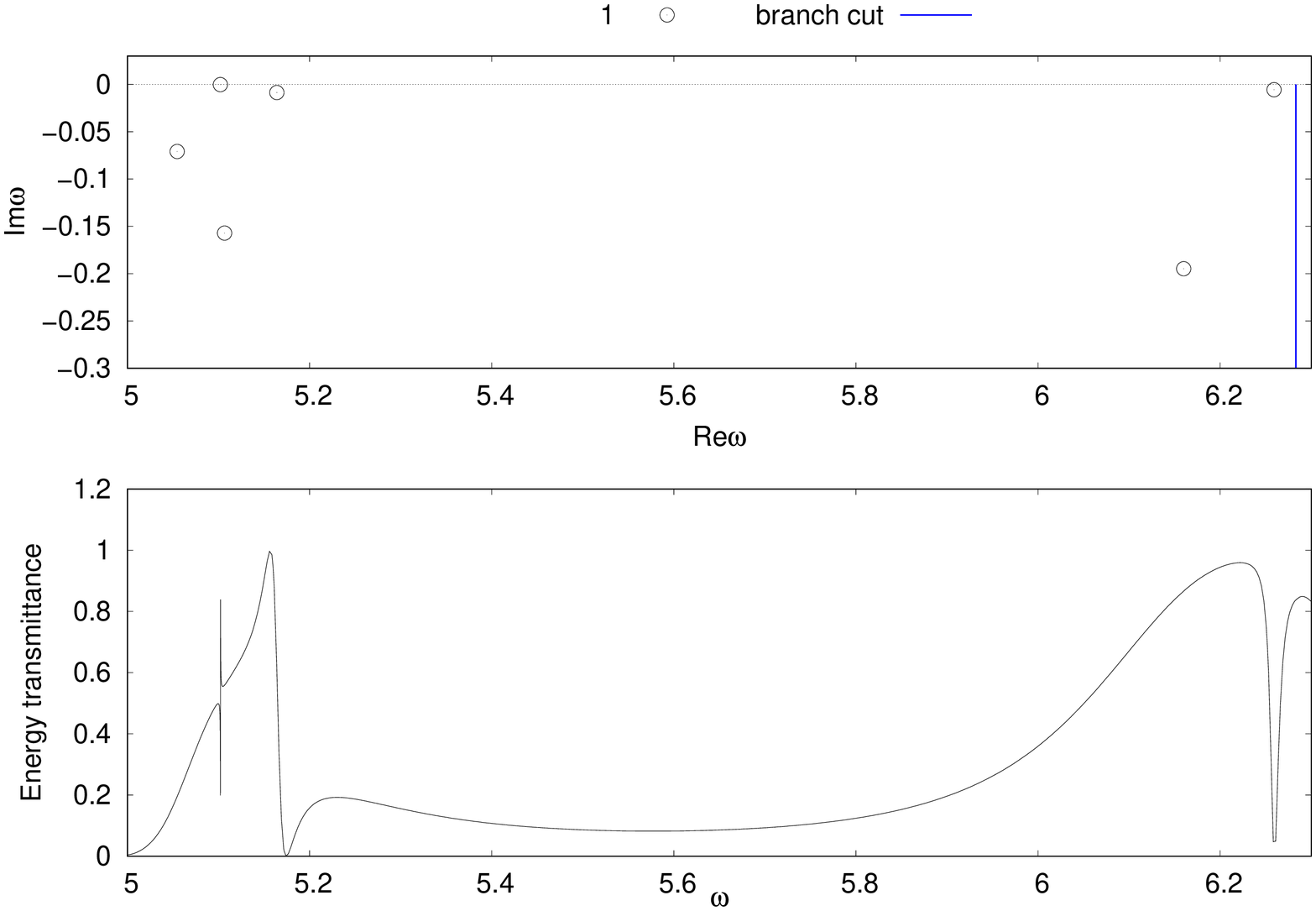}}
     \caption{ Eigenvalues with symmetric modes (upper) and energy transmittance for incident $e^{ik_1x_2}$ (lower) for $\omega\in(5,6.3)$. (Blowup of Fig.\,\ref{fig:eig_even}.)}
     \label{fig:eig_blowup}
   \end{center}
 \end{figure}

Fig.\,\ref{fig:eig} shows all eigenvalues obtained with the 4 methods, where open (solid) symbols stand for the true (fictitious) eigenvalues.
The paths of integration $\gamma$ in (\ref{mu})  are also shown in Fig.\,\ref{fig:eig}.  
This figure shows that the true eigenvalues obtained with the proposed integral equations (methods 1 and 3) agree well with those obtained with the standard approaches (methods 2 and 4) whose accuracy has been examined extensively in our previous paper \cite{Misawaetal}.
This figure further shows that one can clearly distinguish true and fictitious eigenvalues with the proposed methods, in which fictitious ones have positive imaginary parts, while this is not the case with the standard methods (see, for example, the fictitious eigenvalues whose real parts are close to 6).  We note that one of true eigenvalues is very close to the branch point $\pi$,  whose real part, actually,  is slightly smaller than $\pi$. We have checked that this eigenvalue has a sufficiently large reliability index proposed in \cite{Ikegami_Sakurai}.
For validation, we also plot the exact fictitious eigenvalues obtained as the zeros of (\ref{mie_helmholtz}) in Fig.\,\ref{fig:eig}, which agree with numerical fictitious eigenvalues.
  In the upper graph of Fig.\,\ref{fig:eig_even}, we plot only those  true eigenvalues from Fig.\,\ref{fig:eig} obtained with method 1 which have symmetric (with respect to $x_1$)  eigenmodes.  The lower graph of Fig.\,\ref{fig:eig_even} shows 
the energy transmittance produced by  the symmetric incident wave given by $e^{ik_1x_2}$ for real $\omega$. Fig.\,\ref{fig:eig_blowup} is a blow-up of Fig.\,\ref{fig:eig_even} for $\omega\in (5,6.3)$.  We see that the true eigenvalues with small non-positive imaginary parts are close to the peaks or dips  of the energy transmittance for real frequencies.
We observe similar behaviors with eigenvalues having antisymmetric eigenmodes, but we omit the details here.

\begin{figure}
   \begin{center}
   \begin{tabular}{c}
   \begin{minipage}{0.8\textwidth}
     \subfigure[M\"uller]{\includegraphics[width=\textwidth]{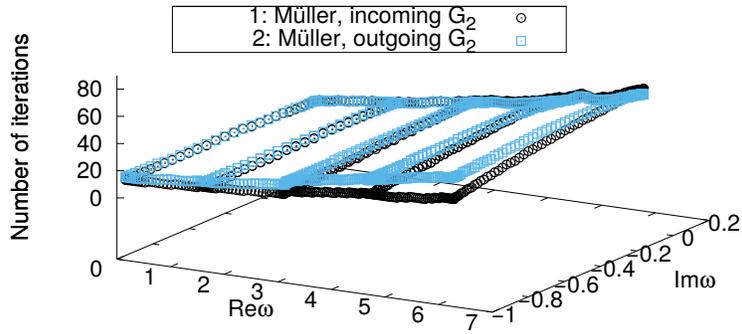}\label{itr_muller}}
   \end{minipage}
   \\
   \begin{minipage}{0.8\textwidth}
     \subfigure[PMCHWT]{\includegraphics[width=\textwidth]{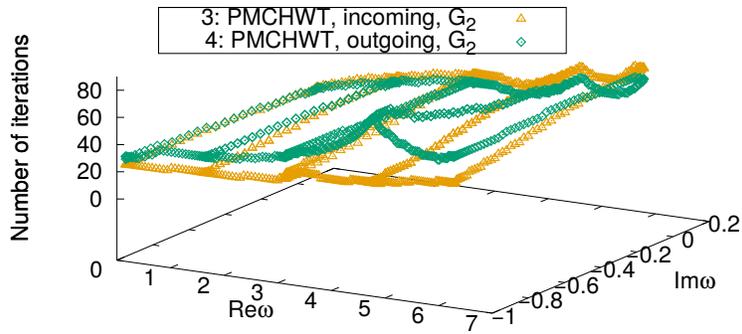}\label{itr_pmchwt}}
   \end{minipage}
   \\
   \begin{minipage}{0.8\textwidth}
     \subfigure[M\"uller, $\Im \omega$ v.s. number of iterations]{\includegraphics[width=\textwidth]{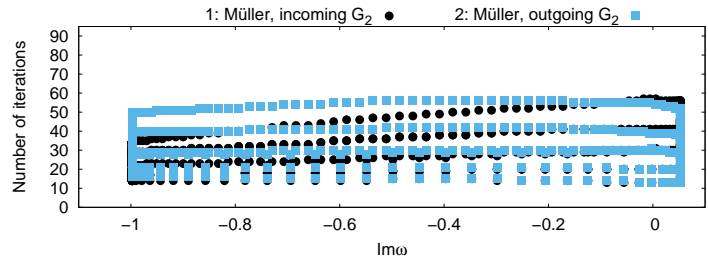}\label{itr_muller_side}}
   \end{minipage}
   \\
  \begin{minipage}{0.8\textwidth}
    \subfigure[PMCHWT, $\Im \omega$ v.s. number of iterations]{\includegraphics[width=\textwidth]{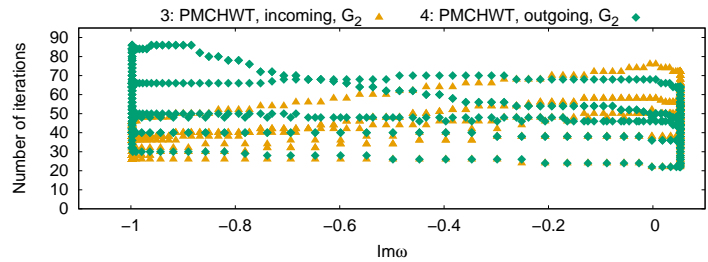}\label{itr_pmchwt_side}}
  \end{minipage}
  \end{tabular}
   \end{center}
   \caption{Number of iterations needed at integration points on the contour $\gamma$.}
  \label{fig:itr}
 \end{figure}

 Fig.\,\ref{fig:itr} shows the number of iterations needed in GMRES  at the integration points on the lower (magenta) paths of integration in Fig.\,\ref{fig:eig} (Note that  the upper paths are used just for obtaining the fictitious eigenvalues of the proposed method,  which are not needed in practice). Figs.\,\ref{fig:itr}\subref{itr_muller} and \subref{itr_pmchwt} show the M\"uller and PMCHWT cases, while Figs.\,\ref{fig:itr}\subref{itr_muller_side} and \subref{itr_pmchwt_side} give the side views ($\Im \omega$ v.s. number of iterations) of Figs.\,\ref{fig:itr}\subref{itr_muller} and \subref{itr_pmchwt}, respectively.
 We observe that the numbers of iterations needed in the proposed methods (methods 1 and 3) are smaller than those of the standard ones when $|\Im \omega|$ is large. 

\clearpage
\subsection{Complex resonances for Maxwell's equations in 3D free space}
We next consider the transmission resonance problem in 3 dimensional free space for the Maxwell equations.  Note that all eigenvalues for the free space transmission problems are leaky and no real eigenvalues exist. However, these problems may have eigenvalues with small imaginary parts, which cause anomalous phenomena as we shall see.

In the first example, we consider a single spherical shell 
\begin{equation*}
\Omega_2=\left\{(x_1,x_2,x_3)\in \mathbb{R}^3\,\,|\,\,0.8< \sqrt{x_1^2+x_2^2+x_3^2}< 1\right\}
\end{equation*}
which encloses $\Omega_3=\left\{(x_1,x_2,x_3)\in \mathbb{R}^3\,\,|\,\,\sqrt{x_1^2+x_2^2+x_3^2}< 0.8\right\}$  and is surrounded by $\Omega_1=\mathbb{R}^3\setminus \overline{\Omega_2\cup \Omega_3}$. We set $\epsilon_1=1, \ \epsilon_2=5, \ \epsilon_3=1$ in $\Omega_1, \  \Omega_2, \ \Omega_3$, respectively, and $\mu=1$ everywhere.
Many true and fictitious eigenvalues appear in this problem, as we shall see.

We now consider the transmission problem for this system. Fictitious eigenvalues for this problem are the eigenvalues for:
\begin{itemize}
  \item a single spherical scatterer $\Omega_3$ having the permitivity $\epsilon=5$ with $\mbox{radius}=0.8$ in the free space having the permitivity $\epsilon=1$
  \item a single spherical scatterer $\Omega_2\cup \Omega_3$ having the permitivity $\epsilon=1$ with $\mbox{radius}=1.0$ in the free space having permitivity $\epsilon=5$.
\end{itemize}
Of course, these fictitious BVPs are with appropriate radiation conditions.

\begin{figure}[]
    \begin{center}
      {\includegraphics[width=\textwidth]{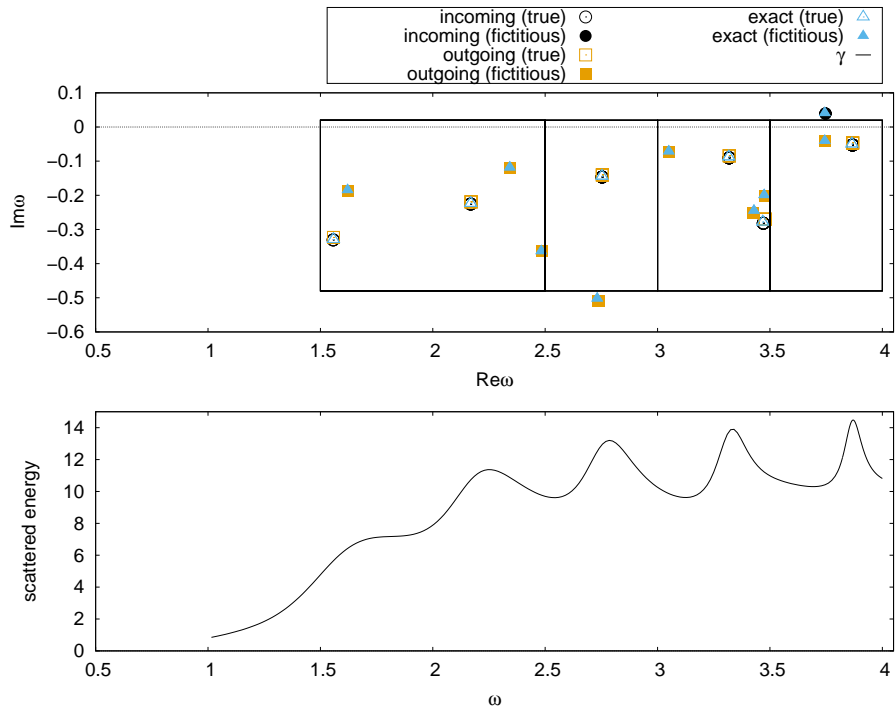}}
      \caption{Upper: eigenvalues, Lower: scattered energy  defined in (\ref{energy}) for incident electric field: $(e^{ik_1 x_3},0,0)$. Solid symbols stand for fictitious eigenvalues.}
      \label{eig_rcs_maxwell}
    \end{center}
  \end{figure}

The upper figure of Fig.\,\ref{eig_rcs_maxwell} shows all eigenvalues obtained with the proposed and standard methods. Here, we used the paths $\gamma$ shown in Fig.\,\ref{eig_rcs_maxwell} and discretized the shell surface with 7300 triangular elements (3920 (3380) elements for $\partial \Omega_1$ ($\partial \Omega_3$) and the total DOF is 87600).  
We set the number of integration points on each side of $\gamma$ to be 16, and the number of random vectors $L$ for SSM to be 20 except for the most left rectangle in the upper figure where $L=10$, respectively. 
  We set larger $L$ and smaller $\gamma$'s  for calculating higher eigenvalues which have larger multiplicities due to geometric symmetry of the problem. One may set a smaller $L$ for geometrically less symmetric scatterers. Note that some eigenvalues outside $\gamma$ are obtained, which happens occasionally in SSM~\cite{Ikegami_Sakurai}

To validate our results, we
note that both true and fictitious eigenvalues for this problem  can be determined easily by means of the Mie-series \cite{Monk}.  These exact eigenvalues, both true ones and fictitious ones for the proposed and standard methods, are plotted with triangular  symbols in the upper figure of Fig.\,\ref{eig_rcs_maxwell}. We see that numerical eigenvalues are obtained correctly and the  proposed method has no fictitious ones in the lower complex plane.

The connection between physical phenomena and true eigenvalues is examined next.
We consider the transmission problem for the same scatterer as above with the incident electric field of $\bm E^\mathrm{inc}=(e^{ik_1 x_3},0,0)$ with real $\omega$  and plot the  scattered energy defined by
\begin{align}
  E^\mathrm{sca}=\int_{\partial \Omega_1} \Re \left( \bm E^\mathrm{sca}\times \bm n\cdot \overline{\bm H^\mathrm{sca}}\right) dS=\int_{\partial \Omega_1} \Re \left( \bm m^\mathrm{sca}\cdot \overline{\bm j^\mathrm{sca}}\times {\bm n}\right) dS
\label{energy}
\end{align}
in the lower figure of Fig.\,\ref{eig_rcs_maxwell} where the superscript ``sca'' stands for the scattered field (i.e., $\bm E^\mathrm{sca}=\bm E-\bm E^\mathrm{inc}$).
We see that  true eigenvalues with imaginary parts smaller than about 0.2 correspond to peaks of the energy.

\begin{figure}
   \begin{center}
   \begin{tabular}{c}
   \begin{minipage}{0.8\textwidth}
     \centering\subfigure[]{\includegraphics[width=\textwidth]{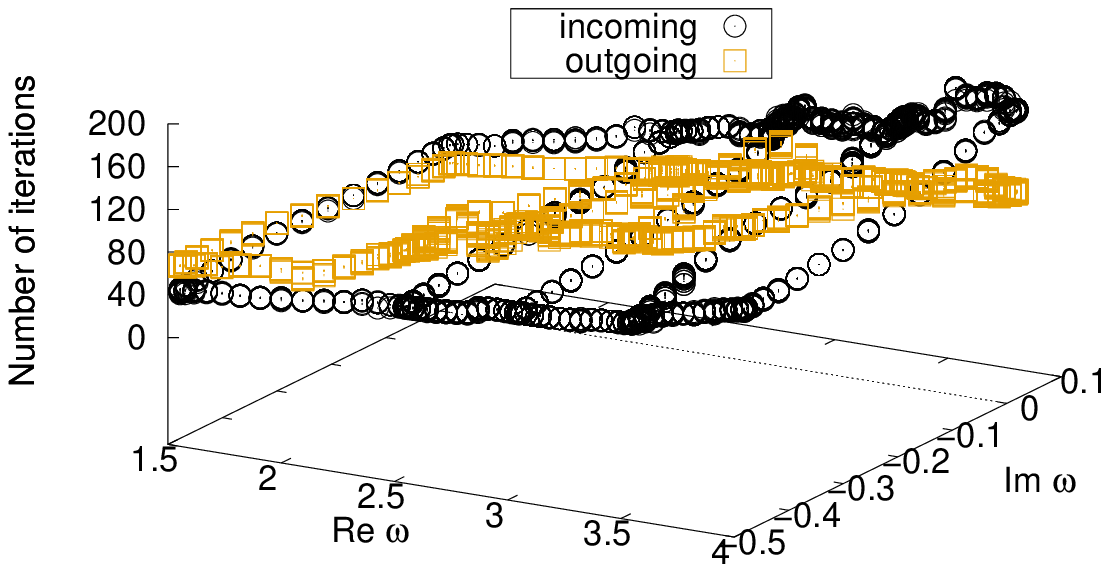}\label{3dplot}}
   \end{minipage}
   \\
   \begin{minipage}{0.8\textwidth}
     \centering\subfigure[]{\includegraphics[width=\textwidth]{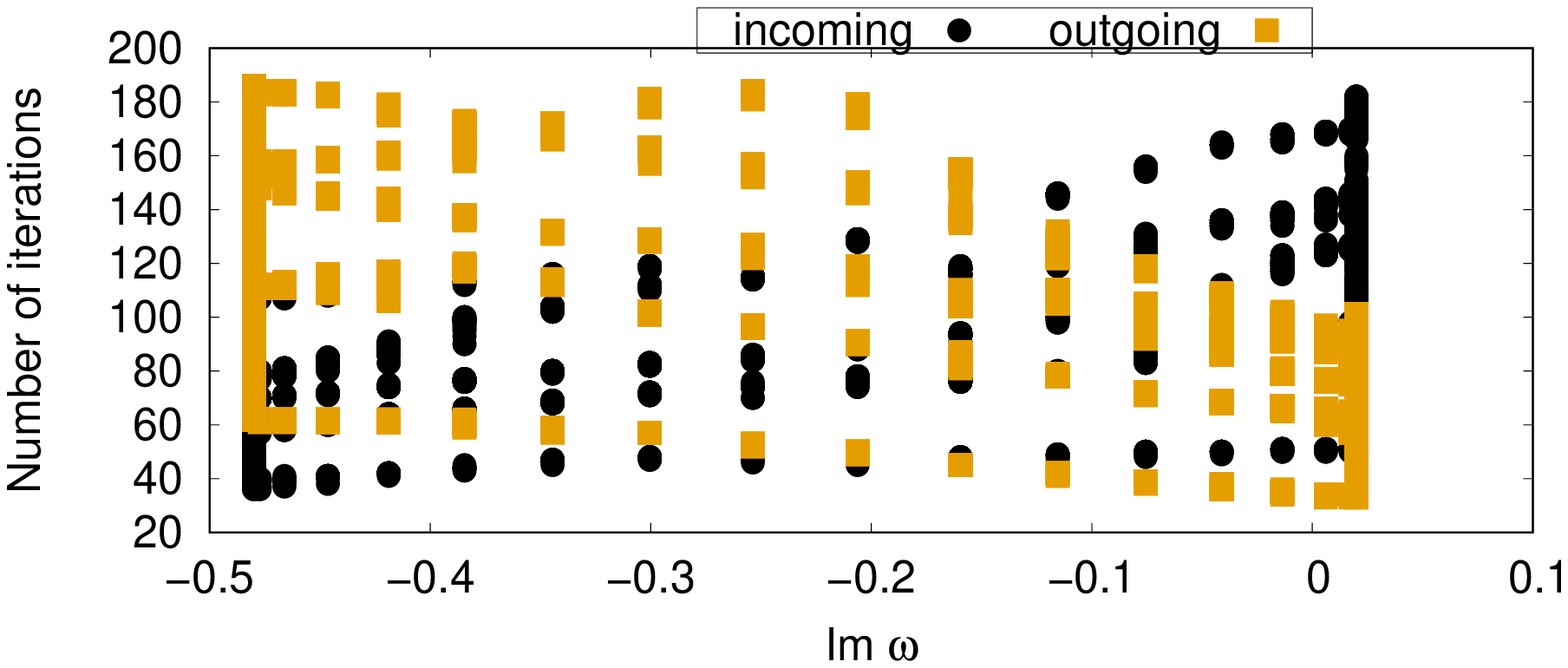}\label{sideview}}
   \end{minipage}
   \end{tabular}
   \caption{Number of iterations needed at integration points on the contour. (spherical shell)}   
  \label{fig:itr_maxwell}
  \end{center}
 \end{figure}

Fig.\,\ref{fig:itr_maxwell}\subref{3dplot} shows the number of iterations needed in  GMRES at the integration points on the paths of integration and Fig.\,\ref{fig:itr_maxwell}\subref{sideview} gives the side view of Fig.\,\ref{fig:itr_maxwell}\subref{3dplot}.  As with the Helmholtz case, we see that the number of iterations needed in the proposed (incoming) method are smaller than those of the standard (outgoing) method when $|\Im \omega|$ is large. However, this is not necessarily the case when $|\Im \omega|$ is small.

The second example is related to eigenvalues for multiple scatterers. We consider two spherical scatterers $\Omega_2$  and $\Omega_3$ ($\Omega_2\cap \Omega_3=\phi$) in the free space $\Omega_1=\mathbb{R}^3\setminus \overline{\Omega_2\cup \Omega_3}$, whose radii are 0.8, 0.4  and centers are $(0,0,0)$, $(1.4,0,0)$, respectively. We set $\epsilon_1=4$, $\epsilon_2=1$ and $\epsilon_3=20$ in $\Omega_1$, $\Omega_2$ and $\Omega_3$, respectively.  The fictitious eigenvalues for this problem  can be obtained easily by means of the Mie-series although the true ones are not very easy to determine.

\begin{figure}[]
     \begin{center}
       {\includegraphics[width=\textwidth]{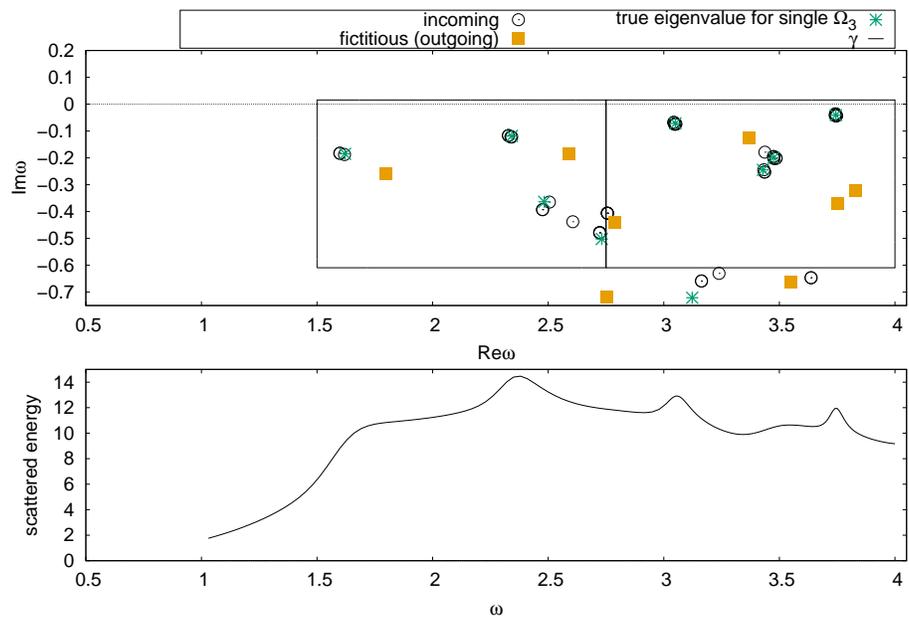}}
       \caption{Upper: eigenvalues, Lower: scattered energy  defined in (\ref{energy}) for incident electric field: $(e^{ik_1 x_3},0,0)$. Note that the asterisks  stand for true eigenvalues for one single scatterer $\Omega_3$.}
       \label{eig_rcs_maxwell_2sph}
     \end{center}
   \end{figure}

The upper figure of Fig.\,\ref{eig_rcs_maxwell_2sph} shows eigenvalues obtained with the proposed (incoming) method. Here, we discretized the surfaces of $\Omega_2$ and $\Omega_3$ with 2880 and 4500 triangular elements, respectively (note that the wavenumber in $\Omega_3$ is higher than that of $\Omega_2$) and the total DOF is 88560. We used 16 (32) integration points on each side of $\gamma$ for the left (right) paths in Fig.\,\ref{eig_rcs_maxwell_2sph} and set $L=20$.
The solid rectangles indicate the (exact) fictitious eigenvalues for this problem which one would obtain with the standard (outgoing) method. This figure clearly shows the usefulness of our method without which it would be very cumbersome, if not impossible,  to separate true eigenvalues from so many fictitious ones. 
Also plotted in the upper figure of Fig.\,\ref{eig_rcs_maxwell_2sph} are the exact true eigenvalues for the single scatterer problem for $\Omega_3$ (i.e., the case with $\epsilon_1=\epsilon_2=4, \ \epsilon_3=20$).  (True eigenvalues for the single scatterer $\Omega_2$ (the case with $\epsilon_1=\epsilon_3=4, \ \epsilon_2=1$) do not exist in the frequency range considered in this figure). We see that true eigenvalues near the real axis for the two spherical scatterer problem are very close to certain eigenvalues for the single scatterer problem, thus indicating that these two scatterer eigenvalues can be interpreted as perturbations of single scatterer eigenvalues.

The lower figure of Fig.\,\ref{eig_rcs_maxwell_2sph} shows the scattered energy in (\ref{energy}) for the incident  electric field given by $\bm E^\mathrm{inc}=(e^{ik_1 x_3},0,0)$ with real $\omega$.  Also in this case, we see that the scattered energy and the eigenvalues with small imaginary parts are related, the latter being close to the peaks of the scattered energy.

\begin{figure}[]
   \begin{center}
     {\includegraphics[width=\textwidth]{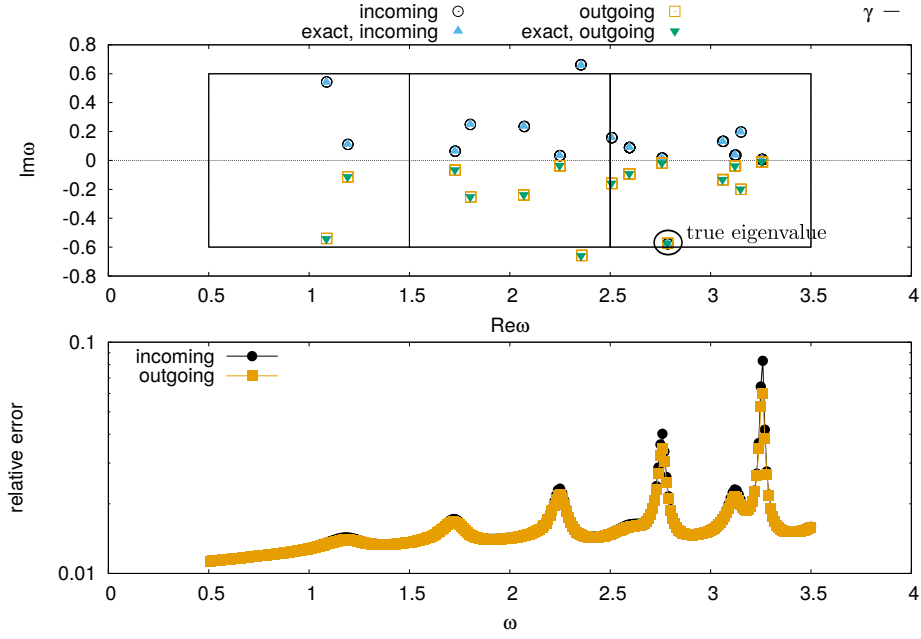}}
     \caption{Upper: eigenvalues, Lower: error  defined in (\ref{err}) of the solution of the transmission problem with incident electric field $(e^{ik_1 x_3},0,0)$.}
     \label{fig:eig_err}
   \end{center}
 \end{figure}

Finally, we show a result which implies that the fictitious eigenvalues may affect the accuracy of BIE solutions in ordinary problems  even when the frequency is real. We consider a single spherical scatterer $\Omega_2$ whose radius is 1.0  and set $\Omega_1=\mathbb{R}^3\setminus \overline{\Omega_2}$. The permitivities in $\Omega_1$ and $\Omega_2$  are $\epsilon_1=6$ and $\epsilon_2=1$, respectively.
We computed  the eigenvalues for this single scatterer problem with the proposed  (incoming) and  standard (outgoing) methods using 5120 triangular elements (61440 DOF) to approximate the surface of the spherical scatterer.
  Also, the paths of integration $\gamma$ for SSM are taken so that they include fictitious eigenvalues with positive imaginary parts just for the purpose of checking. 
We  set the parameter $L$ (see sec.\,\ref{ssmethod}) large ($L=32$) because all eigenvalues for this problem have  large multiplicities.   
For example, the right rectangle in the upper figure of Fig.\,\ref{fig:eig_err} contains as many as 50 eigenvalues. The number of integration points on each side of $\gamma$ is 16.

The upper figure of Fig.\,\ref{fig:eig_err} shows  all the computed eigenvalues,  of which  the one  at $\omega=2.785-0.574i$ (surrounded by a circle) is a true one and others are fictitious.   This figure also includes exact eigenvalues obtained with the Mie-series thus showing that we can determine both  true and fictitious eigenvalues accurately. Note that there exist a few fictitious eigenvalues with very small imaginary parts. We next solved a transmission problem for the same scatterer with the incident electric field given by $\bm E^\mathrm{inc}=(e^{ik_1 x_3},0,0)$ with real $\omega$ and plotted the error relative to the exact solution in the lower figure of Fig.\,\ref{fig:eig_err}, where we defined the error as:
\begin{align}
  \mbox{error}=\frac{1}{2}\left( \frac{\| \bm m-\bm m^\mathrm{Mie}\|}{\|\bm m^\mathrm{Mie}\|}+\frac{\| \bm j-\bm j^\mathrm{Mie}\|}{\|\bm j^\mathrm{Mie}\|}\right)\label{err}.
\end{align}
In (\ref{err}),  $\bm m^\mathrm{Mie}$ and $\bm j^\mathrm{Mie}$ represent the Mie-series solutions, and the norm is the ${\cal L}^2$ norm on the boundary. We see that the error is large near  fictitious eigenvalues with small imaginary parts
regardless of whether we use the proposed integral equation or the standard one.
This result implies that even the M\"uller formulation, which is free of real fictitious eigenvalues,  may  possibly be inaccurate near complex fictitious eigenvalues with small imaginary parts.  Similar observation has been reported by the present authors in \cite{Misawa_Muller}. The same conclusion is quite likely to be true with the PMCHWT formulation which has the same eigenvalues as the M\"uller formulation. 

\section{Concluding remarks}
The results presented in this paper can be summarized as follows:
 \begin{itemize}
   \item We have proposed new boundary integral equations for determining complex eigenvalues of the transmission problems with which one can distinguish true and fictitious eigenvalues easily. We verified that the proposed method could separate the true eigenvalues from the fictitious ones in two dimensional waveguide problems for the Helmholtz equation and three dimensional transmission problems for Maxwell's equations.
 
\item The number of iterations needed in the proposed method for solving BIEs is smaller than that for the standard method when $|\Im \omega|$ is large. However, the situation may be reversed if $|\Im \omega|$ is small. 

 \item Even the proposed  formulation cannot resolve the inaccuracy caused by the fictitious eigenvalues in the ordinary BVPs because it cannot avoid the presence of fictitious eigenvalues close to the real axis.  
\end{itemize}

It is an interesting subject of future work to see if one can find highly accurate BIEs having fictitious eigenvalues with larger imaginary parts than the currently available ones.
The use of fast direct solvers \cite{Martinsson_Rokhlin, Bebendorf} in conjunction with the proposed BIE and SSM is also worth studying because one may sometimes have to take $L$ large, thus having to invert the same matrix repeatedly with many different RHSs. The fast direct solvers may also resolve the problem of increased iteration numbers for small $|\Im \omega|$ mentioned above.
Other future subjects include extension of the proposed method to other problems such as elasticity, etc.

\section*{Acknowledgement}
This work has been supported by JSPS KAKENHI Grant Number 14J03491.

\end{document}